\date{} 
\title{Gaussian free fields for mathematicians}
\author{Scott Sheffield\thanks{Courant Institute.  Partially
supported by NSF grant DMS0403182.}}

\documentclass[12pt,naturalnames]{article}
\usepackage{amsmath}
\usepackage{amssymb}
\usepackage{amsthm}
\usepackage{amsfonts}
\usepackage{graphicx}
\input epsf.tex

\newif\ifhyper\IfFileExists{hyperref.sty}{\hypertrue}{\hyperfalse}

\ifhyper\usepackage{hyperref}\fi

\numberwithin{equation}{section} \numberwithin{figure}{section}

\newtheorem{theorem}{Theorem}
\numberwithin{theorem}{section}

\newtheorem{proposition}[theorem]{Proposition}

\theoremstyle{remark}\newtheorem{definition}[theorem]{Definition}
\theoremstyle{remark}\newtheorem{remark}[theorem]{Remark}

\def\L{\mathcal{L}}

\def\SLEkk#1/{$\mathrm{SLE}_{#1}$} \def\SLEk/{\SLEkk{\kappa}/} \def\SLEtwo/{\SLEkk2/}
\def\SLE/{$\mathrm{SLE}$}
\def\GLEkk#1/{$\mathrm{GLE}_{#1}$} \def\GLEk/{\GLEkk{\kappa}/} \def\GLEtwo/{\GLEkk2/}
\def\GLE/{$\mathrm{GLE}$}
\def\Ito/{It\^o}

\def\R{\mathbb{R}}
          
 \def\qed{\quad \vrule
height7.5pt width4.17pt depth0pt}

\begin{document} \maketitle \begin{abstract} The $d$-dimensional Gaussian free field (GFF),
also called the (Euclidean bosonic) massless free field, is a
$d$-dimensional-time analog of Brownian motion.  Just as Brownian
motion is the limit of the simple random walk (when time and space
are appropriately scaled), the GFF is the limit of many
incrementally varying random functions on $d$-dimensional grids. We
present an overview of the GFF and some of the properties that are
useful in light of recent connections between the GFF and the
Schramm-Loewner evolution.

\end{abstract}

\newpage
\tableofcontents
\newpage

\medbreak {\noindent\bf Acknowledgments.} Many thanks to Oded
Schramm and David Wilson for helping to clarify the definitions and
basic ideas of the text and for helping produce the computer
simulations.  Thanks to Marek Biskup, Yuval Peres, Gabor Pete, Oded
Schramm, Herbert Spohn, Wendelin Werner, and David Wilson for
reading and suggesting improvements to early drafts of this survey
and for recommending additional references.

\section{Introduction}

The $d$-dimensional {\em Gaussian free field} (GFF) is a natural
$d$-dimensional-time analog of Brownian motion.  Like Brownian
motion, it is a simple random object of widespread application and
great intrinsic beauty. It plays an important role in statistical
physics and the theory of random surfaces, particularly in the case
$d=2$. It is also a starting point for many constructions in quantum
field theory \cite{MR1219313, MR707525, MR887102}.

The main purpose of this paper is to provide a mostly self-contained
mathematical introduction to the GFF for readers familiar with basic
probability (Gaussian variables, $\sigma$-algebras, Brownian motion,
etc.), but not necessarily versed in the language of quantum field
theory or conformal field theory.  We will review the classical
continuum constructions (Dirichlet quadratic forms, abstract Wiener
spaces, Gaussian Hilbert spaces, Schwinger functions, chaos
decomposition, etc.) and assemble basic facts about discrete
Gaussian free fields.

Several results from this paper are cited in a recent work by the
author and Schramm, which studies contour lines of the discrete
Gaussian free field and shows that their scaling limits are forms of
the Schramm-Loewner evolution \SLEkk4/ \cite{math.PR/0605337}.  We
also expect these facts to be cited in forthcoming work relating
\SLEk/ to the GFF for other values of $\kappa$.

Although \cite{math.PR/0605337} is a long and technical work, it
contains an elementary twenty-page introduction with many additional
references to the history of the contour line problem and many other
pointers to the physics literature. To avoid duplicating this
effort, we will not discuss \SLE/ at any length here. We also will
not discuss the Virasoro algebra or the use of the GFF in the
Coulomb gas theory (topics discussed at length in several reference
texts, including \cite{MR1052933, MR1424041, MR1694135}), and we
generally make no attempt to survey the physics literature here.
Although this work is primarily a survey, we will also present
without references several simple results (including the natural
coupling of harmonic crystals with the GFF via finite elements and
the coupling of the GFF and Brownian motion via ``field
exploration'') that we have not found articulated in the literature.

\begin{remark} In the physics literature, what we call the GFF is often called the
{\em massless free field} or the {\em Euclidean bosonic massless
free field} --- or else introduced without a title as something like
``the field whose action is the Dirichlet energy'' or ``the Gaussian
field with point covariances given by Green's function.''
\end{remark}

\section{Gaussian free fields}
\subsection{Standard Gaussians}

Consider the space $H_s(D)$ of smooth, real-valued functions on
$\mathbb R^d$ that are supported on a compact subset of a domain $D
\subset \mathbb R^d$ (so that, in particular, their first
derivatives are in $L^2(D)$). This space has a {\em Dirichlet inner
product} defined by $(f_1,f_2)_{\nabla} = \int_D (\nabla f_1 \cdot
\nabla f_2)\,dx$. Denote by $H(D)$ the Hilbert space completion of
$H_s(D)$.  (The space $H(D)$ is in fact a {\em Sobolev space},
sometimes written $\mathbb H_0^1(D)$ or $W^{1,2}_0(D)$
\cite{MR0450957}.) The quantity $(f,f)_{\nabla}$ is called the {\em
Dirichlet energy} of $f$.

Let $g$ be a bijective map from $D$ to another domain $D'$.  If $g$ is a translation or an
orthogonal rotation, then it is not hard to see that $$\int_{D'} \nabla (f_1 \circ g^{-1} ) \cdot
\nabla (f_2 \circ g^{-1} )\,dx = \int_{D} (\nabla f_1 \cdot \nabla f_2)\,dx.$$ If $g(x) = cx$ for a
constant $c$, then $$\int_{D'} \nabla (f_1 \circ g^{-1} ) \cdot \nabla (f_2 \circ g^{-1} )\,dx =
c^{d-2}\int_{D} (\nabla f_1 \cdot \nabla f_2)\,dx.$$

In the special case $d=2$, the equality holds without the $c^{d-2}$
term.  In fact, an elementary change of variables calculation
implies that this equality holds for any conformal map $g$ and any
$f_1, f_2 \in H(D)$. (It is enough to verify this for $f_1, f_2 \in
H_s(D)$.) In other words, the Dirichlet inner product is invariant
under conformal transformations when $d=2$.  (This is one reason
that the GFF is a useful tool in the study of conformally invariant
random two dimensional fractals like \SLE/ \cite{math.PR/0605337}.)

When $D$ is a geometric manifold without boundary (e.g., the unit
torus $\mathbb R^d / \mathbb Z^d$), we define $H_s(D)$ to be the set
of all {\em zero mean} smooth functions on $D$, and again we take
$H(D)$ to be its completion to a Hilbert space with the Dirichlet
inner product.

Note that by integration by parts, $(f_1, f_2)_{\nabla} = (f_1,
-\Delta f_2)$, where $\Delta$ is the Dirichlet Laplacian operator
and $(\cdot, \cdot)$ is the standard inner product for functions on
$D$.  Throughout this paper, we use the notation $(f_1,f_2)_{\nabla}
:= \int_D (\nabla f_1 \cdot \nabla f_2)\,dx$ and $(f_1, f_2) :=
\int_D (f_1 f_2) \, dx$ when the integrals clearly make sense (even
if $f_1$ and $f_2$ do not necessarily belong to $H(D)$ and $L^2(D)$,
respectively).  We also write $\|f\| = (f,f)^{1/2}$ and
$\|f\|_\nabla = (f,f)_\nabla^{1/2}$.

Given any finite-dimensional real vector space $V$ with (positive
definite) inner product $(\cdot, \cdot)$, denote by $\mu_V$ the
probability measure $e^{-(v,v)/2} Z^{-1} dv$, where $dv$ is Lebesgue
measure on $V$ and $Z$ is a normalizing constant.  The following is
well known (and easy to prove) \cite{MR1474726}:

\begin{proposition}\label{p.easygaussian} Let $v$ be a
Lebesgue measurable random variable on $V = \R^d$ with inner product
$(\cdot, \cdot)$ as above. Then the following are equivalent:
\begin{enumerate} \item $v$ has law $\mu_V$.
\item $v$ has the same law as $\sum_{j=1}^d \alpha_j v_j$
where $v_1, \ldots, v_d$ are a deterministic orthonormal basis for
$V$ and the $\alpha_j$ are i.i.d.\ Gaussian random variables with
mean zero and variance one.
\item The characteristic function of $v$ is given by
$$\mathbb E \exp\left( i(v,t)\right) = \exp (- \frac{1}{2}\|t\|^2)$$ for
all $t \in \R^d$.
\item  For each fixed $w \in V$, the inner product
$(v,w)$ is a zero mean Gaussian random variable with variance
$(w,w)$.
\end{enumerate}
\end{proposition}

A random variable satisfying one of the equivalent items in
Proposition \ref{p.easygaussian} is called a {\bf standard Gaussian
random variable on $V$}.  Roughly speaking, the GFF is a standard
Gaussian random variable $h$ on $H(D)$. Because $H(D)$ is infinite
dimensional, some care is required to make this precise. One might
naively try to define $h$ as a random element of $H(D)$ whose
projections onto finite dimensional subspaces of $H(D)$ are standard
Gaussian random variables on those subspaces. However, it is easy
see that this is impossible. (Expanded in terms of an orthonormal
basis, the individual components of $h$ would have to be i.i.d.\
Gaussians --- and hence a.s.\ the sum of their squares would be
infinite, implying $h \not \in H(D)$.)

We will now review two commonly used (and closely related) ways to
define standard Gaussian random variables on infinite dimensional
Hilbert spaces: the abstract Wiener space approach and the Gaussian
Hilbert space approach.


\subsection{Abstract Wiener spaces}
\label{s.abstractwiener}

One way to construct a standard Gaussian random variable $h$ on an
infinite dimensional Hilbert space $H$, proposed by Gross in 1967,
is to define $h$ as a random element not of $H$ but of a larger
Banach space $B$ containing $H$ as a subspace \cite{MR0212152}. To
this end, \cite{MR0212152} defines a norm $| \cdot |$ on $H$ to be
{\bf measurable} if for each $\epsilon
> 0$, there is a finite-dimensional subspace $E_{\epsilon}$ of $H$ for which $$E
\perp E_{\epsilon} \implies \mu_E \left( \{ x \in E: \,\,
|x|>\epsilon \} \right) < \epsilon,$$ where $\mu_E$ is the standard
Gaussian measure on $E$.  In particular, we may cite the following
proposition \cite{MR0212152}.  (Throughout this subsection, $(\cdot,
\cdot)$ and $\| \cdot \|$ denote the inner product and norm of $H$.)

\begin{proposition} \label{p.hilbertschmidt} If $T$ is a {\bf Hilbert Schmidt operator} on $H$ (i.e., $$\sum
\|T f_j\|^2 < \infty$$ for some orthonormal basis $\{f_j\}$ of $H$),
then the norm $\|T \cdot \|$ is measurable.
\end{proposition}

Write $B$ for the Banach space completion of $H$ under the norm $|
\cdot |$, $B'$ for the space of continuous linear functionals on
$B$, and $\mathcal B$ for the smallest $\sigma$-algebra in which the
functionals in $B'$ are measurable. Since each element of $B'$ is a
continuous linear functional on $H$, we may view $B'$ as a subset of
$H$.  Thus $B' \subset H \subset B$. When $b \in B$ and $f \in B'$,
we use the inner product notation $(f,b)$ to denote the value of the
functional $f$ at $b$. (When $f \in H$, this is equal to the inner
product of $f$ and $b$ in $H$.) Given any finite dimensional
subspace $E$ of $B'$ with $H$-orthonormal basis $v_1, \ldots, v_k$,
the map $\phi_E: B \rightarrow E$ given by $\phi_E(b) = \sum (v_j,
b) v_j$ is an extension to $B$ of the orthogonal projection map from
$H$ to $E$. Let $\mu_E$ be the standard Gaussian measure on $E$.
Gross proved the following:

\begin{theorem} If $| \cdot |$ is measurable, then there is a unique probability measure $P$
on $(B, \mathcal B)$ for which $P(\phi_E^{-1} S) = \mu_E(S)$ for each finite dimensional subspace
$E$ of $B'$ and each Lebesgue measurable $S \subset E$.
\end{theorem}

By Proposition \ref{p.easygaussian}, we can restate this as follows:

\begin{theorem} \label{t.grossthm} If $| \cdot |$ is measurable, then there is a unique probability measure $P$
such that if $h$ is a random variable with probability measure $P$
then for any $f \in B'$, the random variable $(h,f)$ is a
one-dimensional Gaussian of zero mean and variance $(f,f)^2$.
\end{theorem}

The triple $(H, B, P)$ is called an {\bf abstract Wiener space}. The
example that motivated Gross's construction is the standard Wiener
space, in which $H = H\left((0,1)\right)$, endowed with the
Dirichlet inner product, $|\cdot|$ is the supremum norm, and $B$ is
the set of continuous functions on $[0,1]$ that vanish on $\{0,
1\}$.  Using the Hilbert space $H(D)$ (with Dirichlet inner product)
we can now give a definition:

\begin{definition} \label{normGFFdefinition}
Given a measurable norm $|\cdot|$ on $H(D)$ and $B'$, $\mathcal B$,
$B$ as above, the {\bf Gaussian free field determined by norm
$|\cdot|$} is the unique $B$-valued, $\mathcal B$-measurable random
variable $h$ with the property that for every fixed $f \in B'$, the
random variable $(h,f)_{\nabla}$ is a Gaussian of variance
$\|f\|_{\nabla}$.  Equivalently, $h = \sum \alpha_j f_j$, where
$\alpha_j$ are i.i.d. Gaussians of unit variance and zero mean and
the $f_j$ are elements of $B$ which form an orthonormal basis for
$H(D)$ --- and the sum is defined within the space $B$. (It is not
hard to see that the partial sums $\sum_{j=1}^m \alpha_j f_j$
converge almost surely in $B$ \cite{MR0212152}.)
\end{definition}

\begin{remark} We can analogously define the {\em complex Gaussian free field
determined by norm $|\cdot|$} by replacing $H$, $B'$, and $B$ with
their complex analogs and writing $h = h_1 + ih_2$, where the $h_1$
and $h_2$ are independent real Gaussian free fields.
\end{remark}

\subsection{Choosing a measurable norm}

We now construct one natural family of measurable norms on $H(D)$
using the eigenvalues of the Laplacian. Suppose that $\{e_j\}$ are
eigenvectors of the Dirichlet Laplacian on $D$ which form an
orthonormal basis of $L^2(D)$ endowed with the usual inner product
and have negative eigenvalues $\{\lambda_j\}$ (ordered to be
non-increasing in $j$). Then an orthonormal basis for $H(D)$ is
given by $f_j = (-\lambda_j)^{-1/2}e_j$, since integration by parts
implies $(e_j,e_k)_{\nabla} = (e_j, -\Delta e_k) = 0$ whenever $j
\not = k$ and $(f_j, f_j)_{\nabla} = ((-\lambda_j)^{-1/2} e_j,
(-\lambda_j)^{1/2} e_j)_\nabla = 1$. (This choice of the $f_j$ is
{\it not} invariant under conformal transformations of $D$ when $d =
2$.)

The reader may recall that by Weyl's formula, if $D \subset \R^d$ is
bounded, then $j^{2/d}/(-\lambda_j)$ tends to a constant as $j \to
\infty$.  (References and much more precise estimates on the growth
of $\lambda_j$ are given in \cite{MR2140257}.)  We define powers of
the negative Dirichlet Laplacian by writing, for each $a \in \R$,
$$(-\Delta)^a \sum \beta_j e_j := \sum (-\lambda_j)^a \beta_j e_j,$$
a definition which makes sense even when $a$ is not an integer.  We
then formally define $\L_a(D) := (-\Delta)^a L^2(D)$ to be the set
of sums of the form $\sum \beta_j e_j$ for which $\sum \beta_j
(-\lambda_j)^{-a} e_j \in L^2(D)$. (When $a<0$, this sum $\sum
\beta_j e_j$ may not converge in $L^2(D)$, but since
$(-\lambda_j)^{-a}$ is polynomial in $j$, it always converges in the
space of distributions on $D$; see Remark \ref{r.distribution}
below.)

Since integration by parts gives $$(f, g)_{\nabla} = \left( f,
(-\Delta) g \right) = \left( (-\Delta)^{1/2} f, (-\Delta)^{1/2} g
\right),$$ the map $(-\Delta)^{-1/2}$ gives a Hilbert space
isomorphism from $L^2(D)$ (with the $L^2$ inner product) to $H(D)$
(with the Dirichlet inner product). Thus we may write $H(D)
=\L_{-1/2}(D)$.

Similarly, for any $a \in \R$, we may view $\L_a(D)$ as a Hilbert
space whose inner product $(\cdot, \cdot)_a$ is the pullback of the
$L^2$ inner product, i.e., $(f,g)_a = ((-\Delta)^{-a} f,
(-\Delta)^{-a} g)$.  We abbreviate $\|f\|_a := ((-\Delta)^{-a} f,
(-\Delta)^{-a} f)^{1/2}$ for the corresponding norm.  An equivalent
way to define $\L_a(D)$ is as the Hilbert space closure of $H_s(D)$
under this norm.

\begin{proposition} \label{p.lapmeas}
Suppose $D$ is a bounded domain in $\R^d$. Then we have the
following:

\begin{enumerate}
\item In the space of formal sums $\sum \beta_j e_j$ (or the space of
distributions) we have $\L_a(D) \subset \L_b(D)$ whenever $a < b$.
\item $\|\cdot \|_b$ is a measurable norm on $\L_a(D)$ (where the latter has inner product
$(\cdot,\cdot)_a$) whenever $a < b - d/4$.
\item When $f \in \L_{-a}(D)$, the functional $g \to (f,g)$ is continuous on $\L_a(D)$.
\end{enumerate}
\end{proposition}

\begin{proof}
The first item is immediate since $\sum
(-\lambda_j)^{-2a}|\beta_j|^2 < \infty$ implies $\sum
(-\lambda_j)^{-2b}|\beta_j|^2 < \infty$.  To prove the second, we
first write $\|f\|_b = \|T_{b-a} f\|_a$ where $T_c := (-
\Delta)^{-c}$. Let $\{g_j\}$ be an orthonormal basis for $\L_b$
under the inner product $(\cdot, \cdot)_b$.  Then $\|f\|_b$ is a
Hilbert Schmidt operator (and hence measurable by Proposition
\ref{p.hilbertschmidt}) provided that $\sum \| T_a g_j \|^2_b = \sum
(-\lambda_j)^{2a-2b} < \infty$. Weyl's formula implies that this
holds provided that $\sum j^{2(2a-2b)/d} < \infty$, which in turn
holds whenever $2(2a-2b)/d < - 1$, i.e., $a < b - d/4$. The final
statement in the proposition is trivial since $$(f,g) = (
(-\Delta)^a f, (-\Delta)^a g)_a = (f, (-\Delta)^{2a} g)_a,$$ and
$(-\Delta)^{2a}g \in \L_a$.\qed
\end{proof}

Proposition \ref{p.lapmeas} implies that, although we cannot
construct the GFF $h$ as a random element of $H(D)$, we can
construct $h$ as a random element of $B=\L_b(D)$, provided $b
> \frac{d-2}{4}$, using the abstract Wiener space definition given
in Section \ref{s.abstractwiener}.

In particular, when $d=1$, we may take $b=0$ and define $h$ as a
random element of $L^2(D)$.  When $d=2$, we cannot define $h$ as a
random element of $L^2(D)$ (indeed, from the power series expansion,
we expect the $L^2$ norm of $h$ to be almost surely infinite), but
we can define $h$ as a random element of $B=\L_b(D)$ for any $b >
0$. In this case, we may view $(h,\cdot)$ as a random continuous
linear functional on $\L_{-b}(D) \subset L^2(D)$ for any $b>0$. In
general, $\rho \to (h,\rho)$ is a random continuous linear
functional on $\L_{-b}(D)$ whenever $b
> \frac{d-2}{4}$.

\begin{remark} \label{r.distribution} Sometimes it is convenient to restrict attention to smooth,
compactly supported test functions $\rho$.  Following the usual
definition, we say that $h$ is a {\bf distribution} if $(h, \cdot)$
is well defined as a functional on the space $H_s(D)$ of smooth
compactly supported functions and this functional is continuous with
respect to the topology of uniform convergence of all derivatives.
If $\rho \in H_s(D)$, then $(-\Delta)^a \rho \in H_s(D) \subset
L^2(D)$ for each positive integer $a$, and it follows that $\rho \in
\L_a(D)$ for all $a$.  If $h \in \L_b(D)$ for some $b$, then $(h,
\cdot)$ is a continuous functional on $\L_{c}(D)$ for any negative
integer $c < -b$.  This implies that the restriction of $(h, \cdot)$
to $H_s(D)$ is continuous in the topology of uniform convergence of
all derivatives (since uniform convergence of all derivatives in
particular implies convergence in $\L_c(D)$), so $h$ is also a
distribution. Many texts (e.g., \cite{MR887102}) simply define the
GFF to be the random distribution determined in this way.  Since
$H_s(D)$ is dense in each of the larger spaces $L_b(D)$, we don't
lose any information by restricting $(h, \cdot)$ to smooth
functions, since there is a unique way to extend $(h, \cdot)$ to a
continuous function on the larger space.
\end{remark}

\begin{remark}
Let $\phi$ be rotationally symmetric smooth positive bump function
on $\R^d$ whose integral is $1$ and which vanishes outside of the
unit ball in $\R^d$.  Let $f_{r,z}(x) = r^{-d} \phi((x-z)/r)$.  This
function is the density of a probability measure on the ball
$B_r(z)$ of radius $r>0$ centered at $z \in \R^d$. Let $D_r$ be the
set of pairs $(r,z)$ for which $r>0$ and $B_r(z) \subset D$. Then
the map $(r,z) \to f_{r,z}$ is continuous from $D_r$ to $\L_b(D)$
for any $b \in \R$. Hence if $h$ is an instance of the GFF defined
by one of the norms discussed above, then $\psi((r,z)) = (h,
f_{r,z})$ is a random continuous function from $D_r$ to $\R$.
Similar arguments show that all the derivatives of $\psi$ are
continuous almost surely. Since the span of the $f_{r,z}$ is dense,
$h$ is almost surely determined by the random smooth function
$\psi$.
\end{remark}

\subsection{Gaussian Hilbert spaces}

The definition of the GFF in terms of abstract Wiener spaces has an
aesthetic and practical drawback in that the choice of measurable
norm is somewhat arbitrary, and it does not yield a description of
the random variable $(h,f)_\nabla$ for general $f \in H(D)$.  In
this section we give a way to make $(h,\cdot)_\nabla$ well defined
as a random variable for each $f \in H(D)$---accepting, of course,
the fact that $f \to (h, f)_\nabla$ cannot be defined as a
continuous functional.

Consider the probability space $(\Omega, \mathcal F, \mu)$ where
$\Omega$ is the set of real sequences $\alpha = \{\alpha_j\}, j \geq
1$, $\mathcal F$ is the smallest $\sigma$-algebra in which the
coordinate projections $\alpha \rightarrow \alpha_j$ are measurable,
and $\mu$ is the probability measure in which the $\alpha_j$ are
i.i.d.\ Gaussian variables of unit variance and zero mean.

In the previous section, we defined the Gaussian free field (GFF) to
be the formal sum $h = \sum_{j=1}^{\infty} \alpha_j f_j$ (which
converges in a larger space $B$), where the $f_j$ are an ordered
orthonormal basis for $H(D)$ and the $\alpha_j$ are i.i.d.\
Gaussians. Now, for any fixed $f \in H(D) = \sum \beta_j f_j$, the
inner product $(h,f)_{\nabla}$ is a random variable that {\em can}
be almost surely well defined as the limit of the partial sums
$\sum_{j=1}^k \beta_j \alpha_j$. (It is important here that we fix
the order of summation in advance, since the sequence $\beta_j
\alpha_j$ is not necessarily a.s.\ absolutely summable.)

Now we have a formal definition:


\begin{definition} \label{GFFdef}  The {\bf Gaussian free field derived from the ordered
orthonormal basis $\{f_j\}$} is the indexed collection $\mathcal
G(D)$ of random variables $(h, f)_{\nabla}$ described above.
\end{definition}

A more abstract definition, which does not specifically reference a
basis or an ordering, is as follows. First, we take the following
definition from \cite{MR1474726}:

\begin{definition} A {\bf Gaussian linear space} is a real linear space of random variables,
defined on an arbitrary probability space $(\Omega, \mathcal F,
\mu)$, such that each variable in the space is a centered (i.e.,
mean zero) Gaussian.  A {\bf Gaussian Hilbert space} is a Gaussian
linear space which is complete, i.e., a closed subspace of
$L^2_{\mathbb R}(\Omega, \mathcal F, \mu)$, consisting of centered
Gaussian variables, which inherits the standard $L^2_{\mathbb
R}(\Omega, \mathcal F, \mu)$ inner product: $(X,Y) = \int XY d\mu$.
We also assume that $\mathcal F$ is the smallest $\sigma$-algebra in
which these random variables are measurable.
\end{definition}

Note that if $X_1, \ldots, X_n$ are any real
random variables with the property that all linear combinations of
the $X_j$ are centered Gaussians, then the joint law of the $X_j$ is
completely determined by the covariances $\text{Cov}[X_j, X_k] =
\mathbb E(X_j X_k)$, and it is a linear transformation of the
standard normal distribution.  A similar statement holds for
infinite collections of random variables \cite{MR1474726}. Then we
have:

\begin{definition} \label{aGFFdef} A {\bf Gaussian free field} is any Gaussian Hilbert space
$\mathcal G(D)$ of random variables denoted by
``$(h,f)_{\nabla}$''---one variable for each $f \in H(D)$---that
inherits the Dirichlet inner product structure of $H(D)$, i.e.,
$$\mathbb E [(h,a)_{\nabla} (h, b)_{\nabla}] = (a,b)_{\nabla}.$$  In
other words, the map from $f$ to the random variable $(h,
f)_{\nabla}$ is an inner product preserving map from $H(D)$ to
$\mathcal G(D)$. \end{definition}

By the identity $(a,b) = \frac{1}{2}[(a+b,a+b) - (a,a) - (b,b)]$, this map is inner product
preserving if and only if it is norm-preserving --- i.e., the variance of $(h, f)_{\nabla}$ is
$(f,f)_{\nabla}$ for each $f \in H(D)$ --- and linear.  Thus we have the following:

\begin{proposition} An $H(D)$-indexed linear space of random variables denoted $(h,f)_{\nabla}$ is a
Gaussian free field if and only if the map from $f\in H(D)$ to the random variable $(h,f)_{\nabla}$
is linear and each $(h,f)_{\nabla}$ is a centered Gaussian with variance $(f,f)_{\nabla}$.
\end{proposition}

Throughout the remainder of this text, we will adopt the Gaussian
Hilbert space approach and view the variables $(h,f)_\nabla$ as
being well defined for all $f \in H(D)$.  Equivalently, we view
$(h,\rho)$ as being well defined for all $\rho \in (-\Delta) H(D) =
\L_{1/2}$.

\begin{remark} When $\rho_1$ and $\rho_2$ are in $H_s(D)$, the covariance of
$(h, \rho_1)$ and $(h, \rho_2)$ can
be written as $(-\Delta^{-1} \rho_1, -\Delta^{-1} \rho_2)_{\nabla} =
(-\Delta^{-1} \rho_1, \rho_2)$.  Since $-\Delta^{-1} \rho$ can be
written using the Green's function kernel as $$[-\Delta^{-1}\rho](x)
= \int_{D} G(x,y) \rho(y) dy,$$ we may also write:

$$\text{Cov}[(h, \rho_1), (h,\rho_2)] = \int_{D \times D} \rho_1(x) \rho_2(y) G(x,y) dxdy$$

When $\rho_1 = \rho_2 = \rho$, the above expression has an
interpretation in electrostatics as the {\bf energy of assembly} of
an electric charge density $\rho$ (grounded at $\partial D$), and
$\Delta^{-1} \rho$ is the {\bf electrostatic potential} of that
density.  The energy of assembly of a density of charge is the
amount of energy required to move charge into that configuration
starting from a zero-energy configuration (in which the potential is
everywhere zero).  Thus the Laplacian $\mathfrak{p} = (-\Delta) h$
is, at least intuitively, a random electrostatic charge distribution
in which the probability of $\mathfrak{p}$ is proportional to $$\exp
( - \text{energy of assembly of $\mathfrak{p}$}).$$ (See
\cite{math.PR/0605337} for more references relevant to this
interpretation.)

\end{remark}

\subsection{Simple examples} \label{torussubsection}
Let $D$ be the unit torus $\mathbb R^d / \mathbb Z^d$.  As before
$H_s(D)$ is the set of smooth functions on $D$ with zero mean and
$H(D)$ is the Hilbert space closure of $H_s(D)$ using the Dirichlet
inner product.  An orthonormal basis for the complex version of
$H(D)$ is given by eigenvectors of the Laplacian, which have the
form $f_k(x) = \frac{1}{2 \pi |k|} e^{2 \pi ix \cdot k}$, for $k \in
\mathbb Z^d \backslash \{0 \}$. Thus, the complex GFF on $D$ is a
random distribution whose Fourier transform consists of i.i.d.\
complex Gaussians times $(2 \pi |k|)^{-1}$.

If $d \geq 2$, then for any fixed $x \in D$ and any fixed ordering
of the $k$'s, the partial sums of $\sum_{j=1}^k \alpha_j f_j(x)$
diverge almost surely, since the variance of the partial sums are
given by $(2 \pi)^{-2} \sum |k|^{-2}$, and this sum diverges when $d
\geq 2$.

When $d = 1$, the limit $h$ can be defined a.s.\ at any given $x$
and is a complex Gaussian.  Since the above sum converges to $(2
\pi)^{-2} 2 \zeta(2) = 1/12$, the real and imaginary components of
$h(x)$ each have variance $1/12$.  In fact, it is not hard to see
that the difference $h(x) -h(0)$ can be written as $(h,
f^x)_{\nabla} = (h, \delta_x - \delta_0)$ where $f^x = - \Delta^{-1}
(\delta_x - \delta_0)$ is continuous and linear on $(0,x)$ and
$(x,1)$.  By computing dot products of $(f^x, f^y)_{\nabla}$, the
reader may verify that $h$ has the same law as a multiple of the
Brownian bridge on the circle, normalized by adding a constant so
that it has zero mean. A similar argument shows that the
one-dimensional GFF on an interval is a multiple of the Brownian
bridge on that interval, and an even simpler argument shows that the
one-dimensional GFF on $(0, \infty)$ is a Brownian motion.  In the
latter case, we may take $[-\Delta^{-1} \delta_x](y) = \min \{x,y
\}$ and $G(x,y) = \min \{x, y\}$.  The variance of $h(x) =
(h,\delta_x)$ is $G(x,x) = x$.

\begin{remark}The definition of the Dirichlet inner product, and hence the
Gaussian free field, has an obvious analog any manifold on which the
Dirichlet energy can be defined.   In particular, since the
Dirichlet inner product is conformally invariant when $d=2$, the
Dirichlet energy has a canonical definition for Riemann surfaces.
There is also a ``free boundary'' version in which we replace
$H_s(D)$ by the set of all smooth, mean zero functions on $D$ with
first derivatives in $L^2(D)$.
\end{remark}

\begin{remark}The GFF has a natural dynamic analog in which each $(h_t,
f)_{\nabla}$ is the Ornstein-Uhlenbeck process with zero mean whose
stationary distribution has variance $\|f\|_{\nabla}^2$.  Thus,
instead of taking $(h, f_j)_{\nabla}$ to be an i.i.d.\ sequence of
random variables, we take $(h_t,f_j)_{\nabla}$ to be an i.i.d.\
sequence of Ornstein-Uhlenbeck processes parameterized by $t$.
\end{remark}

\subsection{Field averages and the Markov property}

If $ - \Delta a = \rho$ is constant on an open subset $D' \subset D$ and equal to zero outside of
$D'$ (i.e., $a$ is harmonic outside of $D'$), then we can think of $(h,a)_{\nabla} = (h, \rho)$ as
describing (up to a constant multiple) the mean value of $h$ on $D'$.  We will retain that
interpretation when $h$ is chosen from the Gaussian free field---i.e., we think of $h$ as
fluctuating so rapidly that it is not necessarily even well-defined as a function, but the
``average value of $h$ on $D'$'' is well-defined.

Since Hilbert spaces are self dual, if $\rho$ is any probability
measure on $D$ for which $f \rightarrow \rho f:= \int f d \rho$ is a
continuous linear functional on $H(D)$ (which is the case if and
only if $\sum [\rho f_j]^2 < \infty$), then there is an $f$ for
which $\rho g = (f,g)_{\nabla}$ for all $g \in H(D)$, and we have
$\rho = -\Delta f \in \Delta H(D)$.

For example, if $d=2$ and $\rho$ is the uniform measure on a line
segment $L$ in the interior of $D$, then $\rho h$ is well-defined.
In this case, the reader may check that $f \rightarrow \rho f$ is
continuous on $H_s(D)$.  (The sums $\sum [\rho f_j]^2 < \infty$ can
be computed explicitly when $D$ is a rectangle; continuity then
follows for domains that are subsets of that rectangle.)

Another important observation is that if $H_1$ and $H_2$ are any
closed orthogonal subspaces of $H(D)$, then $(h,\cdot)_{\nabla}$
restricted to these two subspaces is independent.  To be precise,
denote by $\mathcal F_{H_j}$ the smallest $\sigma$-algebra in which
$h \rightarrow (h,f)_{\nabla}$ is a measurable function for each $f
\in H_j$.  Then it is clear that $\mathcal F_{H_1}$ and $\mathcal
F_{H_2}$ together generate $\mathcal F$, and moreover, $\mu$ is
independent on these two subalgebras.

For example, given an open subset $U$ of $D$, we can write $H_U(D)$ for the closure of the set of
smooth functions that are supported in a compact subset of $U$.  If $a \in H_U(D)$ and $b$ is
harmonic in $U$, then integration by parts implies $(a,b)_{\nabla} = (a, -\Delta b) = 0$.  Thus
$H_U(D)$ is orthogonal to the closed subspace $H^\perp_{U}(D)$ of functions that are harmonic on
$U$.

\begin{theorem}
The spaces $H_U(D)$ and $H^\perp_U(D)$ span $H(D)$. \end{theorem}

\begin{proof} To see this it is enough to show that if $f \in H_s(D)$, then $f$ can
be written as $a + b$, with $a \in H_U(D), b \in H^\perp_U(D)$.  Roughly speaking, we would like to
set $b$ to be the unique continuous function which is equal to $f$ outside of $U$ and harmonic
inside of $U$, and then write $a = f-b$.  But in some cases---e.g., if $U$ is the complement of a
discrete set of points---there is no $b$ with this property.  We will give a slightly modified
definition of $b$ and show $b \in H_U(D)$ and $f-b \in H^\perp_U(D)$.

Let $b_{\delta}(x)$ be the expected value of $f$ at the point at which a Brownian motion started at
$x$ first exits the set $U_{\delta}$ of points of distance more than $\delta$ from the complement
of $U$. Then $a_{\delta}(x) = f - b_{\delta}(x)$ is supported on a compact subset of $U$ and is
clearly in $H_U(D)$.  Since $H_U(D)$ is the closure of the union of the $H_{U_{\delta}}$, the
$a_{\delta}$'s---which are projections onto the increasing (as $\delta \rightarrow 0$) subspaces
$H_{U_{\delta}}$---converge to some function $a \in H_U(D)$.  The $b_{\delta}$ thus must converge
to some $b$, and $b \in H_{U}^\perp$ (since the limit of harmonic functions is harmonic), and $f =
a + b$. \qed \end{proof}

For short, we will write $\mathcal F_U = \mathcal F_{H_U}$ and
$\mathcal F^\perp_U = \mathcal F_{H^{\perp}_U}$.  The
$\sigma$-algebra $\mathcal F^\perp_U$ is one in which random
variables of the form $(h,f)_{\nabla} = (h, -\Delta f)$ are
measurable whenever $\Delta f$ vanishes on $U$. Intuitively, it
allows us to measure the ``values'' of $h$ outside of $U$. On the
other hand, $\mathcal F_U$ allows us to measure the ``values'' of
$h$ inside of $U$ modulo the harmonic functions on $U$.  The
independence of the GFF on $\mathcal F_U$ and $\mathcal F^\perp_U$
can be interpreted as saying that {\em given} the values of $h$
outside of $U$, the distribution of the values of $h$ in $U$ is a
harmonic extension of the values of $h$ on the boundary of $U$ {\em
plus} an independent GFF on $U$. This property of the GFF is called
a {\bf Markov property}. It holds, in particular, if $d=2$ and the
complement of $U$ is a simple path in $D$; in this case, the
$\mathcal F^\perp_U$-measurable functions measure the values of $h$
along that path (or at least the average values of $h$ along
subintervals of that path).

If $U$ is closed and $x \in D \backslash U$, then it is not hard to
see that the projection $f_{x, U}$ of $-\Delta^{-1} \delta_x$ onto
$H_U$ has finite Dirichlet energy and that its Laplacian is
supported on the boundary of $D \backslash U$.  Although $h(x) =
(h,\delta_x) = (h,-\Delta^{-1} \delta_x)_{\nabla} $ is not a
well-defined random variable, we may still intuitively interpret
$(h, f_{x,U})_{\nabla}$ as the ``expected'' value of $h(x)$ given
the values of $h$ in $U$.  The reader may check that the function
$(h,f_{x,U})_{\nabla}$ is almost surely harmonic in $D \backslash
U$. (More precisely, since the event ``$(h,f_{x,U})_{\nabla}$ is
harmonic'' is not in our $\sigma$-algebra, one shows that the
function $(h,f_{x,U})_{\nabla}$ defined on dyadic rational points of
$D \backslash U$ almost surely extends continuously to a harmonic
function in all of $D \backslash U$; this is the same way one proves
continuity of Brownian motion in, e.g., Chapter 7 of
\cite{MR1609153}.) We interpret this function as the harmonic
extension to $D \backslash U$ of the ``values'' of $h$ on the
boundary of $U$.

\subsection{Field exploration: Brownian motion and the GFF} \label{exploresection}

In this subsection, we describe a simple way of using a space-filling curve to give a linear
correspondence between the GFF and Brownian motion.  Roughly speaking, we ``explore'' the field $h$
along a space-filling curve, and the Brownian motion goes up or down depending on whether the
values we encounter are greater than or less than what we expect. Then each of the random variables
$(h,f)_{\nabla}$ can be viewed as an appropriate stochastic integral of this Brownian motion.
Although analogous constructions hold in higher dimensions, we will assume for simplicity that
$d=2$ and $D$ is a simply connected bounded domain.

First, choose $f_0$ so that $\Delta f_0$ is a negative constant on
$D$. Then let $\gamma:[0,1] \rightarrow D$ be a continuous
space-filling curve.  For each $t$, denote by $\gamma_t$ the compact
set $\gamma([0,t])$, and let $P_t$ be the projection onto the
subspace $H^\perp_{D \backslash \gamma_t}$ of functions harmonic in
$D \backslash \gamma_t$. We also require that $\gamma$ remains
continuous when it is parameterized in such a way that
$\|P_t(f_0)\|_{\nabla}^2 = t$ for all $t \in
[0,\|f_0\|^2_{\nabla}]$.  (This will be the case provided that
$\gamma_s$ is a proper subset of $\gamma_t$ whenever $s < t$.  In
other words, although $\gamma$ may intersect itself, it cannot spend
an entire positive-length interval of time retracing points that
have already been seen.)

By decomposing $f_0$ into its projection onto the complementary subspaces $\mathcal F^\perp_{D
\backslash \gamma_t}$ and $\mathcal F_{D \backslash \gamma_t}$, we easily observe the following:
$$W(t) := \mathbb E ( (h,f_0)_{\nabla} | \mathcal F^\perp_{D \backslash \gamma_t}) = (h,
P_t(f_0))_{\nabla}.$$ Clearly, $W$ is a martingale, and each $W(t) -
W(s)$ is Gaussian with variance $|s - t|$.  Hence, $W$ has the same
law as a Brownian motion (in the smallest $\sigma$-algebra where
each $W(t)$ is measurable).

Now, the reader may easily verify that the linear span of the
functions $P_t(f_0)$, with $t \in [0,1]$, is dense in $H(D)$.  Thus,
given the Brownian motion $W(t)$, it should be possible, almost
surely, to determine $(h,f_j)_{\nabla}$ for each $j$.  To this end,
observe that for any other $f$, the value $$W_f(t) := \mathbb E (
(h,f)_{\nabla} | \mathcal F^\perp_{D \backslash \gamma_t}) = (h,
P_t(f))_{\nabla}$$ is also a martingale, and is Brownian motion when
time is parameterized by $\|P_t(f)||y_{\nabla}^2$.  The question is,
how are $W_f$ and $W$ related?

The answer would be obvious if we had $f = a P_s(f_0)$ for some
fixed constants $a$ and $0<s<1$. In this case, $W_f(t) = a W (\min
\{s,t \})$.  A similar result holds if $f$ is any finite sum of such
functions.  We have now defined $W_f(t)$ for a dense linear space of
functions $f \in H(D)$, and we may choose an orthonormal basis
$\{f_j\}$ for $H(D)$ from among that space.  Given an arbitrary $f =
\sum \alpha_j f_j$ and any fixed $t$, we can take $W_f(t)$ to be the
limit of the partial sums of $\sum \alpha_j W_{f_j}(t)$.

The above discussion gives a {\em linear} correspondence between the
GFF and a Brownian motion. That correspondences of this sort should
exist is not surprising, given that a Brownian motion can be
interpreted as a Gaussian free field on the interval $I = (0,
\|f_0\|^2_{\nabla})$ that is only required to vanish at the left
endpoint (i.e., $H(I)$ is the Dirichlet inner product Hilbert space
completion of the set of smooth functions on $I$ that vanish at
zero). The map that sends the function $g^s(t) = \min \{s,t\}$ in
$H(I)$ to $P_s(f_0) \in H(D)$ extends to a Hilbert space isomorphism
between $H(D)$ and $H(I)$ and the correspondence between the GFF and
Brownian motion is induced by this isomorphism.

\subsection{Circle averages and thick points}
Fix a domain $D \subset \mathbb R^2$ on which the GFF is defined.
For $x \in D$ and $t \in \mathbb R$, write $B_x(t)$ for the mean
value of the GFF on the circle of radius $e^{-t}$ centered at $x$
(provided $t$ is large enough that the disc enclosed by this circle
lies in $D$).  The reader may verify that for each $x \in D$, the
law of $B_x(t)$ is that of a multiple of a Brownian motion.  If $x
\not = y$, then the GFF Markov property implies that the Brownian
motions $B_x(t)$ and $B_y(t)$ grow independently of one another once
$t$ is large enough so that $2e^{-t} < |x-y|$.  Let $$A_x(t) =
\frac{\int_{s = t}^{\infty} B_x(s) e^{-s}ds}{\int_{s = t}^{\infty}
e^{-s}ds} = \int_{s = t}^{\infty} B_x(s) e^{t-s}ds$$ be the mean
value of the GFF on the disc of radius $e^{-t}$.

One can generate various fractal subsets of $D$ by considering the
set of points $x$ for which $B_x(t)$ or $A_x(t)$ are in some sense
highly atypical.  For example, one might consider the set of $x$ for
which $\limsup |B_x(t)| \leq C$ for a constant $C$, or more
generally the set of $x$ for which $\limsup |B_x(t) - f(t)| \leq C$
for a some function $f$.

In \cite{HuPeres} the authors fix a constant $0 \leq a \leq 2$ and
define a {\bf thick point} as a point $x \in D$ for which $\lim
A_x(t)/t = \sqrt{a}$.  They prove that the Hausdorff dimension of
the set of thick points of a GFF is almost surely equal to $2 - a$.


\section{General results for Gaussian Hilbert spaces}
By replacing the Dirichlet inner product with a different bilinear
form, it is possible to construct different types of Gaussian
Hilbert spaces, some of which play important roles in constructive
quantum field theory \cite{MR887102}.  Most of the results in the
previous section depended heavily on the choice of inner product.
In this section, we will back up and make some statements about the
Gaussian free field that are largely independent of this choice.

\subsection{Moments and Schwinger functions} \label{momentsection}

The moments of the variables $(h, \rho)$ may be computed explicitly.
Suppose that $\rho_j = -\Delta f_j$ for $1 \leq j \leq k$.  First,
we know that $\mathbb E \left[(h,\rho)_{\nabla}\right] = 0$ and
$\mathbb E \left[(h,\rho_1)(h,\rho_2)\right] = (f_1, f_2)_{\nabla}$.
We can now cite the following from the first chapter of
\cite{MR1474726}:
\begin{theorem} \label{momenttheorem} We have $$\mathbb E[(h,\rho_1) \cdots (h,\rho_k)]
= \sum_{M} \prod_{j=1}^{k/2} (f_{M_{j,1}},f_{M_{j,2}})_{\nabla}$$
where $M$ ranges over the set of all partitions $M = \{ (M_{j,1},
M_{j,2}) \}$ of $\{1,\ldots,k\}$ into $k/2$ disjoint pairs. In
particular, this value is zero whenever $k$ is odd. \end{theorem}

If each $(f_j, f_k)_{\nabla}$ is positive (which will be the case,
e.g., if each $\rho_j$ is a positive probability density function),
we can interpret the above value as the partition function for a
process that chooses a random perfect matching $M$ of the complete
graph on $1, \ldots, k$ with probability proportional to
$$\prod_{j=1}^{k/2} (-\Delta^{-1} \rho_{M_{j,1}},\rho_{M_{j,2}})$$
i.e., each edge $(i,j)$ is weighted by $(f_i, f_j)_{\nabla}$. (This
is a useful mnemonic, if nothing else.)  These perfect matchings are
simple examples of {\bf Feynman diagrams} (see \cite{MR1474726} or
\cite{MR887102}).  The moment computations above apply to general
Gaussian Hilbert spaces when $(\cdot ,\cdot)_{\nabla}$ is replaced
by the appropriate covariance inner product.

Now, for any $x_1, \ldots, x_k \in D$, one would also like to define
the point moments $S_k(x_1,\ldots, x_k) := \mathbb E \prod_{j=1}^k
h(x_j)$.  In the case of the GFF when $d \geq 2$, these $h(x_j)$ are
not defined as random variables.  If they {\em were} defined as
random variables, then, writing $\rho = -\Delta f$, we would also
expect that $\mathbb E(h, \rho) = \int_D S_1(x)\rho(x) dx$ and more
generally $$\mathbb E \prod_{j=1}^k (h, \rho_j) = \int_{D^k}
S_k(x_1, \ldots, x_k) \prod \rho_j(x_j) dx.$$ It turns out that for
a broad class of Gaussian and non-Gaussian random fields that
includes the GFF, there {\em do} exist functions (or at least
distributions) $S_k :D^k \rightarrow \mathbb R$, called {\bf
Schwinger functions}, for which the latter statement holds, at least
when the $\rho_j$ are smooth (Proposition 6.1.4 of \cite{MR887102}).

In the case of the GFF, the reader may verify that $S_k$ is identically zero when $k$ is odd and
$S_2(x_1, x_2) = G(x_1,x_2)$. For even $k > 2$, $S_k$ can be computed from $S_2$ using the
expansion given in Theorem \ref{momenttheorem}.

\subsection{Wiener decompositions and Wick products}

From Chapter 2 of \cite{MR1474726}, we cite the following:

\begin{theorem} The set of polynomials in the $(h,f_j)_{\nabla}$ is a dense subspace of $L^2(\Omega,
\mathcal F, \mu)$.  \end{theorem}

The space $L^2(\Omega, \mathcal F, \mu)$, endowed with the inner
product $(X,Y) = \mathbb E (XY)$, can be viewed as the closure of
the direct sum of Hilbert spaces $H^{:n:}$, each of which is the
closure of the set of degree $n$ polynomials in $(h,f_j)_{\nabla}$
that are orthogonal to all degree $n-1$ polynomials (in particular,
$H^{:1:} = H(D)$).  The decomposition of an element of $L^2(\Omega,
\mathcal F, \mu)$ into these spaces is sometimes called the {\bf
Wiener chaos decomposition} (Wiener 1938, It\^{o} 1951, Segal 1956;
see Chapter 2 of \cite{MR1474726}). See Chapters 2 and 3 of
\cite{MR1474726} for the explicit form of the projection operators
onto each $H^{:n:}$, as well as a natural orthonormal basis for each
$H^{:n:}$---defined in terms of an orthonormal basis of $H(D)$.

The above implies in particular that if we use the two-dimensional
Gaussian free field to define conformally invariant random sets or
random loop ensembles (e.g., via \SLE/ constructions), then any
$L^2$ function of these random objects (e.g., the number of loops
encircling a given disc) can be expanded in terms of this
orthonormal basis (although in practice this may be difficult to do
explicitly).

Now, given any $\eta_1 \in H^{:m:}$ and $\eta_2 \in H^{:n:}$, the
{\bf Wick product} of $\eta_1$ and $\eta_2$ is the projection of
$\eta_1 \eta_2$ onto $H^{:m+n:}$.  Explicit formulae for the Wick
product are given in \cite{MR1474726}.

\subsection{Other fields}

The ``massive'' free fields (see Chapter 6 of \cite{MR887102}) may
be defined as a collection of Gaussian random variables $(h,\rho)$
with covariances given by $$\text{Cov}[(h,\rho_1), (h,\rho_2)] = (
(-\Delta + m^2)^{-1} \rho_1, \rho_2),$$ where $m$ is a real number
called the {\bf mass}. Here, we can either require the test
functions $\rho$ to be smooth (the random distribution
interpretation) or let them belong to the Hilbert space completion
of the smooth functions under the inner product $((-\Delta+m^2)^{-1}
\rho_1, \rho_2)$ (the Gaussian Hilbert space interpretation).

Equivalently, we may consider random variables $(h, f)^m_{\nabla} : = (h,f)_{\nabla} + m^2(h,f)$
where $f = (-\Delta + m^2)^{-1} \rho$.  We then have $$\text{Cov}[(h,f_1)^m_{\nabla},
(h,f_2)^m_{\nabla}] = (f_1, f_2)^m_{\nabla}.$$  The GFF is the case $m=0$.  Most of the results in
this paper have straightforward analogs in the case $m \not = 0$.

Among the other fields discussed in quantum field theory (see, e.g.,
Chapters 8, 10, and 11 of \cite{MR887102}) are probability measures
that are absolutely continuous with respect to the massive or
massless free fields and have Radon-Nikodym derivatives given by
$e^{-V} Z^{-1}$, where $V$ is a ``potential'' function on $\Omega$
and $Z$ is an appropriate normalizing constant.

In the more interesting ``interacting particle'' settings, however,
this $V$ is undefined or infinite for $\mu$-almost all points in
$\Omega$.  This happens, for example, if $d \geq 2$ and we write
$V(h) = \int_{D} P(h)$ where $P$ is an even polynomial.  In this
case, we can define approximations $V_n$ to $V$ by writing $V_n(h) =
V(h_n)$ where the $h_n$ is a natural approximation to $h$ (e.g.,
$h_n$ could be the partial sum $\sum_{j=1}^{n} \alpha_j f_j$, or it
could be one of the discrete lattice approximations in the
subsequent section).  We then seek to define a field which is the
limit, in some sense, of (appropriately normalized versions of) the
probability measures $e^{-V_n} Z_n^{-1} \mu$.  When it exists, the
limiting measure is in general not absolutely continuous with
respect to $\mu$.  See \cite{MR887102} for a mathematically rigorous
approach to constructing fields, including the so-called {\bf
$P(\phi)$ fields} (in particular, the celebrated $\phi^4$ fields),
in a way that uses the Gaussian fields as a starting point.

\section{Harmonic crystals and discrete approximations of the GFF}
\subsection{Harmonic crystals and random walks}

Let $\Lambda$ be a finite graph with a positive weight function $w$
on its edges.  If $\phi_1$ and $\phi_2$ are functions on $\Lambda$,
we denote their Dirichlet inner product by $$(\phi_1,
\phi_2)_{\nabla} = \sum_{e=(x,y)} w(e)(\phi_1(y) -
\phi_1(x))(\phi_2(y) - \phi_2(x))$$ where the sum is over all edges
$e=(x,y)$ of $\Lambda$.  Now, fix a ``boundary'' $\partial \Lambda$,
which, for now, we can take to be any non-empty subset of the
vertices of $\Lambda$. Then the set $H(\Lambda)$ of real-valued
functions on $\Lambda$ whose values are fixed to be zero (or some
other pre-determined set of boundary values) on $\partial \Lambda$
is a $|\Lambda| - |\partial \Lambda|$ dimensional Hilbert space
under the Dirichlet inner product. This is also a finite dimensional
Gaussian Hilbert space, where the probability density at $\phi$ is
proportional to $e^{-\|\phi\|_{\nabla}^2/2}$.

When $\Lambda$ is a large subset of a $d$-dimensional lattice graph
$L$, and $\partial \Lambda$ is the set of vertices that border
points of $L \backslash \Lambda$, the resulting ``discretized random
surface'' model is commonly called the {\bf discrete Gaussian free
field} (DGFF), the {\bf discrete massless free field}, or the {\bf
harmonic crystal}.

Next, consider a random walk in which each edge $e=(x,y)$ is activated by an independent
exponential clock with intensity $w(e)$, at which point if the position is $x$, it switches to $y$,
and vice versa.

Given $\phi$ and $x \in \Lambda \backslash \partial \Lambda$, write
$Y$ for the expected value of $\phi(x)$ at the first neighbor $x$ of
$y$ hit by a random walk starting at $y$, i.e., $$Y = \frac{\sum_{e
= (x,y)} w(e) \phi(x)}{\sum_{e = (x,y)} w(e)}.$$ When the weight
function $w$ is constant, $Y$ is simply the average value of $\phi$
on the neighbors of $x$; when $w$ is not constant, $Y$ is an
appropriate weighted average of these values.

The reader may verify the following one-point discrete Markov
property: the random variable $\phi(y) - Y$ is independent of the
values of $\phi(x)$ for $x \not = y$ and this random variable has
mean zero and variance $1/\sum_{e=(x,y)} w(e)$.

This property implies that if we fix the values of $\phi$ on $\partial \Lambda$ (where these values
are not necessarily equal to zero), then the expected value of any $\phi(x)$, for $x \in \Lambda$,
will be the expected value of $\phi(x_h)$, where $x_h$ is the first vertex on $\partial \Lambda$
that is hit by a random walk beginning at $x$.

We claim that it also implies that $\text{Cov} [\phi(x),\phi(y)]$ is
equal to the expected amount of time that a particle started at $x$
will spend at $y$ before hitting the boundary.  (This function in
$x$ and $y$ is also called the {\bf discrete Green's function} on
$\Lambda$.)  To see this, first observe that when $y$ is fixed, both
sides are discrete harmonic (with respect to the random walk) in
$\Lambda \backslash \{y\}$.  It is then enough to compute the
discrete Laplacian (with respect to the random walk) at $y$ itself
and observe that it is equal to $1/\sum_{e \ni y} w(e)$ for both
sides: for the right hand side (the random walk interpretation),
this follows from well-known properties of exponential clocks.  For
the left hand side, since $Y$ and $\phi(y) - Y$ are independent and
$\mathbb E (\phi(y) - Y)^2 = 1/\sum_{e=(x,y)} w(e)$, we have
$$\mathbb E \phi(y)(\phi(y) - Y) = \mathbb E (\phi(y) - Y)^2 +
\mathbb E Y (\phi(y) - Y) = 1/\sum_{e=(x,y)} w(e).$$

As an example, if $\Lambda$ is $[-n,n]^d \subset \mathbb Z^d$ and $w=1$, then it is not hard to see
(from well-known properties of random walks) that the variance of $\phi(0)$ is asymptotically
proportional to $n$ if $d=1$, $\log (n)$ if $d = 2$, and a constant if $d \geq 3$ (since the random
walk is transient in the latter case).

See \cite{MR1880237, MR1757956, MR1468313, GiacominTransworld} and
the references therein for these and many more results about
harmonic crystals and various generalizations of the harmonic
crystal.  For an analog of the Green's function interpretation of
variance that applies in the continuum case (known as the {\bf
Dynkin isomorphism theorem}), see \cite{MR693227, MR734803,
MR756768}.

\subsection{Discrete approximations: triangular lattice}

Suppose $d=2$, and let $L$ be the standard triangular lattice (the
dual of the honeycomb lattice). Now, suppose we restrict the GFF on
$D$ to the $\sigma$-algebra $\mathcal F_{H_n}$ where $H_n(D)$ is the
set of continuous functions that are affine on the each of triangles
of $\frac{1}{n} L$ and that vanish on the boundary of $D$. Since
$H_n(D)$ is a finite dimensional Hilbert space, it is self-dual, and
hence a sample from the GFF determines an element $h_n$ of $H_n(D)$,
with probability proportional to $\exp[-\|h_n\|_{\nabla}^2]$.
Observe also that $h_n$ is determined by its values on the vertices
of the triangular mesh, and that $\|h_n\|_{\nabla}^2$ is equal to
$\frac{\sqrt 3}{6} \sum |h_n(j) - h_n(i)|^2 + |h_n(k)-h_n(i)|^2 +
|h_n(k) - h_n(j)|^2$ where the sum is over all triangles $(i,j,k)$
in the mesh. (The area of each triangle is $\sqrt 3/(4n^2)$ and the
gradient squared is $\frac{2}{3} n^2 [|h_n(j) - h_n(i)|^2 +
|h_n(k)-h_n(i)|^2 + |h_n(k) - h_n(j)|^2]$.) Since each interior edge
of $D$ is contained in two triangles, this is also equal to
$3^{-1/2} \sum |h_n(j) - h_n(i)|^2$, where the sum is taken over
nearest neighbor pairs $(i,j)$.

In other words, $h_n$ is distributed like $3^{1/4}$ times the
harmonic crystal (with unit weights) on the set of vertices of
$\frac{1}{n} L \cap D$, where the boundary vertices are precisely
those that lie on a triangle which is not completely contained in
$D$.

It is not hard to see that the union of the spaces $H_n(D)$ is dense in $H(D)$, and that for any
function $f \in H(D)$, $f$ minus its projection onto $H_n(D)$ will tend to zero in $n$.  Thus, the
$h_n$ are approximations of $h$ in the sense that for any fixed $f \in H(D)$, we have $(h_n,
f)_{\nabla} \rightarrow (h,f)_{\nabla}$ almost surely.  For any fixed $n$, it easy to define
contour lines of the continuous function $h_n$.


\subsection{Discrete approximations: other lattices}

Now, again, suppose that $d=2$ but that we replace the standard lattice with any doubly periodic
triangular lattice $L'$ (i.e., a doubly periodic planar graph in which all faces are triangles).
Once again, we can restrict the GFF to the functions that are linear on the triangles of
$\frac{1}{n}L'$.  In this case, the heights at the corners turn out to have the same law as a
harmonic crystal on this triangular lattice graph in which the weight $w(e)$ corresponding to a
given edge $e$ is given by $w(e) = [\cot(\theta_1) + \cot(\theta_2)]/2$, where $\theta_1$ and
$\theta_2$ are the angles opposite $e$ in the two triangles that are incident to $e$. (Note, of
course, that the weights are unchanged by constant rescaling.)

That is (as the reader may verify), the discrete Dirichlet energy $$(h_n,h_n)_{\nabla} = \sum
w((x,y)) [h_n(y) - h_n(x)]^2,$$ where $h_n$ is viewed as a function on the graph, is the same as
the continuous Dirichlet energy $(h_n,h_n)_{\nabla}$, where $h_n$ is extended piecewise linearly to
all of $D$.

In particular, we can get the discrete GFF on a square grid by dividing each square into a pair of
triangles (either direction).  In this case, the diagonal edges have $w(e) = 0$, since both
$\theta_1$ and $\theta_2$ are right angles.  The vertical and horizontal edges have $w(e)=1$, since
$\theta_1$ and $\theta_2$ are equal to $45$ degrees in this case.  (Note that $w(e)$ can be
negative if one or both of the $\theta_i$ exceeds 90 degrees.  Earlier, we assumed that $w$ was
positive; however, since the quadratic form corresponding to such a $w$ is still positive definite,
the definition of the harmonic crystal still makes sense in this setting.)

Similar results hold when $d > 2$ if we replace triangles with
$d$-dimensional simplices.

\subsection{Computer simulations of harmonic crystals}

One fortunate feature of the discrete Gaussian free fields is the
ease with which they can be simulated on computers.  Consider first
the special case that $\Lambda$ is an $m \times n$ torus grid graph.
In this case, the discrete exponential functions form an orthogonal
basis of the set of complex, mean-zero functions on $\Lambda$, so
the elements in a discrete Fourier transform of a discrete GFF on
$\Lambda$ are independently distributed Gaussians with easily
computed variances. For example, the Mathematica code shown below
generates and plots an instance of the discrete GFF on an $m \times
n$ torus.

{\bf \begin{verbatim}ListPlot3D[ Re[Fourier[Table[ \end{verbatim}
\begin{verbatim}   (InverseErf[2 Random[]-1]+I InverseErf[2 Random[]-1]) * \end{verbatim}
\begin{verbatim} If[j+k==2,0,1/Sqrt[(Sin[(j-1)*Pi/m]^2+ Sin[(k-1)* Pi/n]^2)]],\end{verbatim}
\begin{verbatim} {j,m}, {k,n}]]]] \end{verbatim}}

The code generates a complex GFF and then plots its real part.  To
parse the code, note that the second line simply produces a complex
Gaussian random variable in Mathematica.  For each $1 \leq j \leq
m$, $1 \leq k \leq n$, the third line gives one over the gradient
norm of the discrete exponential $(x,y) \rightarrow \eta^{(j-1) x}
\zeta^{(k-1) y}$ (where $\eta$ and $\zeta$ are $m$th and $n$th roots
of unity, respectively) on the $m \times n$ torus (unless $j=k=1$,
in which case it gives zero). The Fourier transform of the
corresponding matrix is the Gaussian free field on the torus. Now
suppose $\Lambda$ is a simply connected induced subgraph of the
torus with values $h_0$ assigned to its boundary; to sample a GFF on
$\Lambda$ with these boundary values, one may first sample an
instance $h$ of the GFF on the torus and then replace $h$ with
$h+\tilde h$ where $\tilde h$ is the discrete harmonic interpolation
of the function $h_0-h$ (defined on $\partial \Lambda$) to all of
$\Lambda$.  The most time-consuming part of this algorithm is
computing the discrete harmonic interpolation (but it is not hard to
compute an approximate interpolation).

\begin{figure}\label{dgff} \epsfbox[-50 40 150 150]{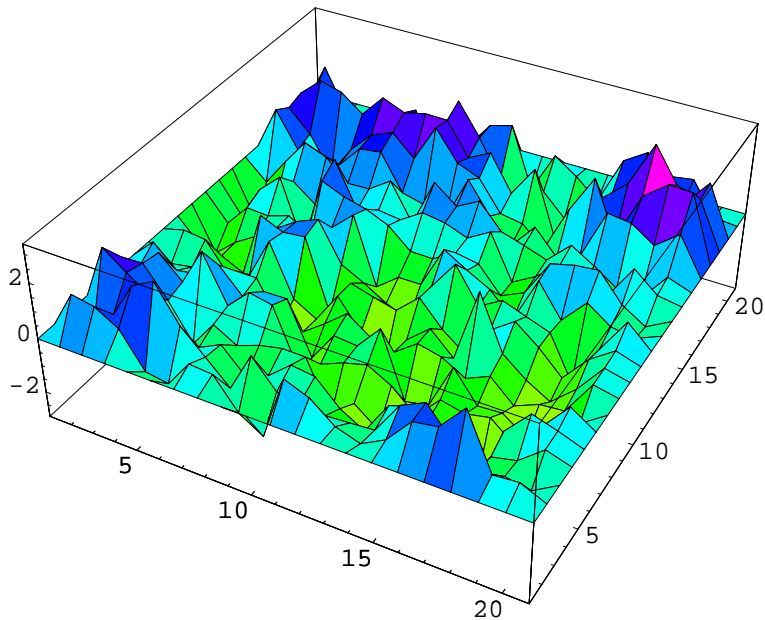} \caption {Discrete Gaussian free field
on 20 by 20 grid with zero boundary conditions.}
\end{figure}

\subsection{Central limit theorems for random surfaces} Kenyon in \cite{MR1872739} proved that random domino
tiling height functions (in regions with certain kinds of boundary
conditions) converge to the GFF as the mesh size gets finer. (To be
precise, \cite{MR1872739} studies random discrete height functions
$h_{\epsilon}$---for which the lattice spacing is $\epsilon$---and
shows that for smooth density functions $\rho$, $(h_{\epsilon},
\rho)$ converges in law to $(h, \rho) = (h,-\Delta^{-1}
\rho)_{\nabla}$.) Also, \cite{MR1461951,MR1872740} give a similar
Gaussian free field convergence result for a class of discretized
random surfaces known as {\bf Ginzburg-Landau $\nabla \phi$ random
surfaces} or {\bf anharmonic crystals} and \cite{MR1872740} shows
further that certain time-varying versions of these processes
converge to the dynamic GFF.

See \cite{SheffieldAsterisque} for more references on random surface
models (assuming both discrete and continuous height values) with
convex nearest neighbor potential functions and discussion of the
division of the so-called {\bf gradient Gibbs measures} into {\bf
smooth} phases (in which the variance of the height difference
between two points is bounded, independently of the distance between
those points) and {\bf rough} phases (in which the variance of these
height differences tends to infinity as the points get further
apart). An important open question is whether every two-dimensional
rough phase has a scaling limit given by a linear transformation of
the GFF.

\bibliographystyle{halpha}
\addcontentsline{toc}{section}{Bibliography}
\bibliography{GFFsurvey}

\bigskip

\filbreak
\begingroup
\small
\parindent=0pt

\bigskip
\vtop{
\hsize=2.3in
Courant Institute\\
New York University\\
New York, NY 10012, USA\\ \\
{\tt {sheff@math.nyu.edu}} } \endgroup \filbreak  \end{document}